\documentclass{amsart}
\newtheorem{theorem}{Theorem}[section]
\newtheorem{lemma}[theorem]{Lemma}

\newtheorem{corollary}[theorem]{Corollary}

\theoremstyle{definition}

\newtheorem{example}[theorem]{Example}

\theoremstyle{remark}
\newtheorem{remark}[theorem]{Remark}

\numberwithin{equation}{section}

\begin{document}
\title{Splitting theorems in presence of an irrotational vector field}
\author{Manuel Guti\'errez\and Benjam\'{\i}n Olea}

\address{Departamento de \'Algebra,
Geometr\'{\i}a y Topolog\'{\i}a. Universidad de M\'alaga, Campus
Teatinos, 29071- M\'alaga, Spain.} \email{M.G: {\ttfamily
mgl@agt.cie.uma.es}\\ B.O: {\ttfamily benji@agt.cie.uma.es}}


\maketitle

\begin{abstract}
New splitting theorems in a semi-Riemannian manifold which admits
an irrotational vector field (not necessarily a gradient) with
some suitable properties are obtained. According to the extras
hypothesis assumed on the vector field, we can get twisted, warped
or direct decompositions. Some applications to Lorentzian manifold
are shown and also $\mathbf{S}^{1}\times L$ type decomposition is
treated.

\end{abstract}

\section{Introduction}
\label{intro}

Warped products are a generalization of direct products, giving
sophisticated examples of semi-Riemannian manifolds from simpler
ones. They are manageable for computations and sufficiently rich
to have a great geometrical and physical interest. The standard
spacetime models of the universe and the simplest models of
neighborhoods of star and black holes are warped products,
therefore, it is of interest to know when a Lorentzian manifold
can be decomposed as a warped product. In this paper, we give
decomposition theorems without assuming simply connectedness nor
the existence of a gradient, obtaining warped and twisted
decomposition. Given two semi-Riemannian manifolds $(B,g_{B})$,
$(L,g_{L})$ and a function $f\in C^{\infty }(B)$, the warped
product $M=B\times _{f}L$ is the product manifold furnished with
the metric $g=g_{B}+f^{2}g_{L}\;$\cite{Oneill}$.$ When $f$ is a
$C^{\infty }$ function on $B\times L$, it is called a twisted
product.

The classical theorems that ensures the metric decomposition of a
manifold as a direct product are the De Rham and De Rham-Wu
decomposition theorem \cite{DeRham}, \cite{Wu}. They were
generalized by Ponge and Reckziegel in \cite{Ponge}, where more
general decompositions, such as twisted and warped products, were
obtained. In all these paper the manifolds are simply connected.

More recent advances, in which a non necessarily simply connected
manifold is decomposed as a product, assume the existence of a
function without
critical points \cite{Fischer}, \cite{Garcia1}, \cite{Garcia2}, \cite{Innami}%
, \cite{Kanai}, \cite{Sakai}. It is a great simplification because
it ensures that the integral curves of the gradient meet the level
hypersurfaces of the function for only one value of its parameter.
This
allows us to construct explicitly a diffeomorphism between the manifold and $%
\mathbf{R}\times L,$ where $L$ is a level hypersurface. Some
additional properties of the gradient permit to get different
types of metric
decompositions. The fact that the function has not critical points exclude $%
\mathbf{S}^{1}\times L$ decompositions$,$ which are not frequent
in the literature.

There are other results in which it is assumed the existence of a
vector field which is not necessarily a gradient. In the paper
\cite{Garcia3} it is
obtained a metric decomposition of a manifold as a direct product $\mathbf{R}%
\times L$ assuming some conditions on a timelike vector field and
its orthogonal distribution. On the other hand, a diffeomorphic
decomposition can be given in a chronological manifold furnished
with a special vector field \cite{Harris}.

Sometimes, the decomposition theorems are stated as singularity
versus splitting theorems: if a manifold is not a global product,
it must be incomplete \cite{Garcia1}, \cite{Garcia2},
\cite{Garcia3}, \cite{Walchap2}.

The decomposition process of a manifold has two stages:
diffeomorphic and metric. One of the standard hypothesis to obtain
a diffeomorphic decomposition is simply connectedness, which
obviously is not a necessary condition.

Once we have a diffeomorphic decomposition $B\times L,$ we can
obtain the metric decomposition assuming some geometrical
properties on the canonical foliation of the product $B\times L$
\cite{Ponge}.

An usual technique used in the literature to split
diffeomorphically a manifold is to construct a diffeomorphism
using the flow of a suitable vector field. Although the
construction of the diffeomorphism is the same in all cases, each
theorem is proved in a different way depending on the hypotesis
assumed on the vector field. In section two it is shown that the
flow of an unitary vector field with an additional property,
present in most splitting theorem in the literature, induces a
local diffeomorphism which is onto. This gives us a common basis
to obtain decomposition theorems. The difficult part is to check
the injectivity, which is equivalent to ensure that each integral
curve of the vector field intersects the leaves of the orthogonal
distribution only in one point.

In section three a general decomposition lemma is presented, which
is the basic tool to obtain the splitting theorems.

In section four some decomposition theorems for irrotational
vector field with compact leaves are given, and in section five
they are applied to Lorentzian manifolds. In section six the
$\mathbf{S}^{1}\times L$ type decomposition is treated.

All manifolds considered in this paper are assumed to be
connected. We follow the sign convention for curvature of
\cite{Beem}, $\left. R_{XY}Z=\nabla _{X}\nabla _{Y}Z-\nabla
_{Y}\nabla _{X}Z-\nabla _{[X,Y]}Z\right. $ and we write $Ric(X)$
for the quadratic form associated
with the Ricci curvature tensor$.$ Given $f:M\rightarrow \mathbf{R}$ a $%
C^{\infty }$ function, we call $grad\,f$ the gradient of $f,$
$H^{f}$ its Hessian and $\triangle f=div\,grad\,f$ its laplacian.

\section{Preliminaries on the flow of an unitary vector field}
\label{sec:1}

Let $(M,g)$ be a semi-Riemannian manifold and $U$ a vector field
on $M$ with never null norm. The vector field $U$ has integrable
orthogonal distribution if and only if it is orthogonally
irrotational, i.e., $g(\nabla _{X}U,Y)=g(X,\nabla _{Y}U)$ for all
$X,Y\in U^{\perp }.$ In this situation, we call $L_{p}$ the leaf
through $p,$ $E$ its unitary, $\lambda $ the function such that
$U=\lambda E$, $\varepsilon =g(E,E)$ and $\Phi $ the flow of $E.$
Usually, the metric that we put on the leaf is the induced metric.
The vector field $U$ is irrotational if $g(\nabla
_{X}U,Y)=g(X,\nabla _{Y}U)$ for all vector fields $X,Y$ on $M$. We
say that $U$ is pregeodesic if its unitary is geodesic$,$ or
equivalently, if $\nabla _{U}U$ is proportional to $U.$

It is useful to know when the flow of an unitary and orthogonally
irrotational vector field takes leaves into leaves, because it
facilitates the construction of a diffeomorphism using the flow
restricted to an orthogonal leaf.

\begin{lemma}
\label{LemaHojasAhojas}Let $M$ be a semi-Riemannian manifold and
$E$ an unitary, orthogonally irrotational and complete vector
field. Then $\Phi _{t}$ satisfies $\left. \Phi _{t}(L_{p})\subset
L_{\Phi _{t}(p)}\right. $ for all $t\in \mathbf{R}$ and $p\in M$
if and only if $\nabla _{E}E=0.$
\end{lemma}
\begin{proof}
 Suppose that $\Phi _{t}$ takes leaves into leaves.
Then, if $v\in E^{\perp }$ it follows that $\left. g(E_{_{\Phi
_{t}(p)}},(\Phi _{t})_{*p}(v))=0\right. $ for
all $t\in \mathbf{R},$ i.e., $(\Phi _{t})^{*}(g)(E_{p},v)=0$ for all $t\in \mathbf{%
R}.$ Then
\begin{equation*}
(L_{E}g)_{p}(E_{p},v)=0.
\end{equation*}

But given a vector field $A\in E^{\perp }$,
\begin{equation*}
L_{E}g(E,A)=g(\nabla _{E}E,A),
\end{equation*}
then $g(\nabla _{E}E,A)=0$ for all $A\in E^{\perp }$ and being $E$ unitary, $%
\nabla _{E}E=0.$

Now suppose $\nabla _{E}E=0$. This implies $E$ is irrotational and
so the metrically equivalent one-form $w$ is closed. Then

\begin{equation*}
L_{E}w=d\circ i_{E}w+i_{E}\circ dw=0,
\end{equation*}
so $\Phi _{t}^{*}w=w.$ Therefore $\Phi _{t}$ takes leaves into
leaves for all $t\in \mathbf{R}.$
\end{proof}

\begin{remark}
\label{Piezas de Hojas a Hojas}We suppose that $E$ is a complete
vector field for convenience. If we do not assume it, we should
say that the flow takes any connected open set of a leaf into a
leaf. Note also that being $E$ unitary and orthogonally
irrotational, it is geodesic if and only if it is irrotational.
\end{remark}

Now, let $U$ be an orthogonally irrotational vector field with
never null norm in a semi-Riemann manifold. If $E$ is complete and
geodesic, since it is orthogonally irrotational too, we can apply
lemma \ref{LemaHojasAhojas}. Take $\Phi $ the flow of $E$ and $L$
an orthogonal leaf. We construct
\begin{eqnarray*}
\Phi &:&\mathbf{R}\times L\rightarrow M \\
&&\left. (t,p)\rightarrow \Phi _{t}(p).\right.
\end{eqnarray*}

Since $(\Phi _{t})_{*p}(E_{p})=E_{_{\Phi _{t}(p)}}$ and $\Phi
_{t}$ takes leaves into leaves, $\Phi $ is a local diffeomorphism
which preserves the foliations and indentifies $E$ with
$\frac{\partial }{\partial t}.$

\begin{lemma}
\label{LemaFlujoSobre}Let $M$ be a semi-Riemannian manifold and
$E$ an unitary, irrotational and complete vector field. Then the
local diffeomorphism $\Phi :\mathbf{R}\times L\rightarrow M$ is
onto.
\end{lemma}
\begin{proof}
We show that $M=\cup _{t\in \mathbf{R}}\Phi _{t}(L).$ It is
sufficient to verify that $\cup _{t\in \mathbf{R}}\Phi _{t}(L)$ is
an open and closed set. It is an open set because we know that
$\Phi $ is a local diffeomorphism. If we
take $x\notin \cup _{t\in \mathbf{R}}\Phi _{t}(L)$, then $x\in \cup _{t\in \mathbf{%
R}}\Phi _{t}(L_{x})\subset (\cup _{t\in \mathbf{R}}\Phi _{t}(L))^{c},$ but $%
\cup _{t\in \mathbf{R}}\Phi _{t}(L_{x})$ is also an open set, so
$\cup _{t\in \mathbf{R}}\Phi _{t}(L)$ is a closed set.
\end{proof}

If we can ensure that $\Phi :\mathbf{R}\times L\rightarrow M$ is
also injective we would have a diffeomorphic decomposition of $M.$
In most of the splitting theorems which we can find in the
literature, the vector field (or its unitary) verifies the
hypothesis of the lemma \ref{LemaFlujoSobre}. The injectivity of
$\Phi $ is equivalent to that the integral curves of $E$ meet the
orthogonal leaves for only one value of its parameter.

\section{Global decomposition lemma}
\label{sec:2}

It is well known that two orthogonally and complementary foliation
give rise to a local diffeomorphic decomposition of the manifold.
Depending on certain geometrical properties of the foliations, we
can get also a metric decomposition.\ Consider $g$ a metric on
$M_{1}\times M_{2}$ such that the canonical foliations are
orthogonal. Take $(p_{0},q_{0})\in M_{1}\times
M_{2} $, $F_{p_{0}}:M_{2}\rightarrow M_{1}\times M_{2}$ given by $%
F_{p_{0}}(q)=(p_{0},q)$ and $F^{q_{0}}:M_{1}\rightarrow
M_{1}\times M_{2}$
given by $F^{q_{0}}(p)=(p,q_{0}).$ Now, we construct the metrics $%
g_{1}=(F^{q_{0}})^{*}(g)$ and $g_{2}=(F_{p_{0}})^{*}(g).$ Then, in
\cite {Ponge} it is proven that

\begin{enumerate}
\item  If both canonical foliations are geodesic, the metric is
the direct product $g_{1}+g_{2}.$

\item  If the first canonical foliation is geodesic and the second
one
spherical, the metric is the warped product $g_{1}+f^{2}g_{2},$ where $%
f(p_{0})=1.$

\item  If the first canonical foliation is geodesic and the second
one
umbilic, the metric is the twisted product $g_{1}+f^{2}g_{2}$, where $%
f(p_{0},q)\equiv 1.$

\item  If the first canonical foliation is geodesic, the metric is
of the form $g_{1}+h_{p},$ where $h_{p}=(F_{p})^{*}(g),$ i.e., for
each $p\in M_{1}, $ $h_{p}$ is a metric tensor on $M_{2}$ (in some
special cases this is called a parametrized product
\cite{Garcia1}).
\end{enumerate}

An orthogonally irrotational vector field with never null norm
gives rise to two orthogonal and complementary foliations. We can
make extra hypothesis about the vector field to obtain geometrical
properties about the foliations and achieve metric decompositions.
We say that a vector field $U$ is orthogonally conformal if there
exists a $\rho \in C^{\infty }(M)$ such that $(L_{U}g)(X,Y)=\rho
g(X,Y)$ for all $X,Y\in U^{\perp }$. The following result codifies
the properties of the foliations in terms of the normalized of the
vector field.

\begin{lemma}
\label{Propiedades de las dos Foliaciones}Let $M$ be a
semi-Riemannian manifold and $U$ an orthogonally irrotational
vector field with never null norm.

\begin{enumerate}
\item  The foliations $U$ and $U^{\perp }$ are totally geodesic if
and only if $E$ is parallel.

\item  The foliation $U$ is totally geodesic and $U^{\perp }$
spherical if and only if $E$ is irrotational, orthogonally
conformal and $grad\,divE$ is proportional to $E$.

\item  The foliation $U$ is totally geodesic and $U^{\perp }$
umbilic if and only if $E$ is irrotational and orthogonally
conformal.
\end{enumerate}
\end{lemma}
\begin{proof}
The case $1$ is trivial and the \textit{if} part of case $2$ can
be found in \cite{Miguel}.
We prove first the third case. Assume that $U$ is totally geodesic and $%
U^{\perp }$ is umbilic. Then $E$ is orthogonally irrotational and
geodesic, and therefore it is irrotational. If we call $II$ the
second fundamental form of $U^{\perp },$ since it is umbilic,
there is $b\in C^{\infty }(M)$ such that $II(X,Y)=g(X,Y)\cdot bE$,
for all $X$, $Y\in E^{\perp }.$ On the other hand,
$II(X,Y)=\varepsilon g(\nabla _{X}Y,E)E=-\varepsilon g(Y,\nabla
_{X}E)E,$ so $\nabla _{X}E=-\varepsilon bX+\alpha (X)E.$ But since
$E$ is
unitary $\alpha (X)=0.$ Therefore $\nabla _{X}E=-\varepsilon bX$ for all $%
X\in E^{\perp }$ and then $E$ is orthogonally conformal. The
converse is easy. If in addition $U^{\perp }$ is spherical then
$b$ is constant through the leaves and so $grad\,divE$ is
proportional to $E$. This prove the \textit{only if} part of case
2$.$
\end{proof}
 If $U$ is an irrotational and conformal vector field,
its unitary verifies the case two of the lemma, and if $U$ is an
irrotational, orthogonally conformal and pregeodesic vector field,
its unitary verifies the case three. In any case, $\lambda $ is
constant through the orthogonal leaves.

It is easy to verify that $U$ is irrotational and conformal if and only if $%
\nabla U=a\cdot id$ for some $a\in C^{\infty }(M).$ In this situation, $%
a=E(\lambda ).$

On the other hand, $U$ is irrotational, orthogonally conformal and
pregeodesic if and only if $\nabla _{X}U=aX+bg(X,E)E,$ where
$a,b\in C^{\infty }(M)$.

The following result is the key to prove the splitting theorems
given is this paper.

\begin{lemma}
\label{ProposFacil}Let $M$ be a semi-Riemannian manifold and $E$
an unitary, irrotational and complete vector field. Take $p\in M$
and suppose that the integral curves of $E$ with initial value on
$L_{p}$ intersect $L_{p}$ at
only one point. Then $M$ is isometric to $\mathbf{R}\times L_{p}$ or $\mathbf{S}%
^{1}\times L_{p}$ with metric $g=\varepsilon dt^{2}+g_{t}$ (a
semi-Riemannian parametric product) where $g_{0}=g\mid _{L_{p}}$,
and $E$ is identified with $\frac{\partial }{\partial t}.$
\end{lemma}
\begin{proof}
Take $\Phi :\mathbf{R}\times M\rightarrow M$ the flow of $E.$ We know that $%
\Phi :\mathbf{R}\times L\rightarrow M$ is a local diffeomorphism
and onto because of lemma \ref{LemaFlujoSobre}. If the integral
curves $\Phi _{t}(q),$ $q\in L_{p},$ meet $L_{p}$ for only one
value of its parameter, then $\Phi $ is injective. If one of them
meet $L_{p}$ again, then all the curves $\Phi _{t}(q),$ $q\in
L_{p},$ meet $L_{p}$ again since $\Phi $ takes leaves into leaves.
We know that the integral curves intersect $L_{p}$ at only one
point, so the curves $\Phi _{t}(q),$ $q\in L_{p},$ must be
periodic. It is easy to verify that they have the same period, let us say $t_{0},$ i.e., $%
\Phi (t_{0},q)=q$ for all $q\in L_{p}.$ Then, we can define a
diffeomorphism
\begin{eqnarray*}
\Psi &:&\mathbf{S}^{1}\times L\rightarrow M \\
&&\left. (e^{it},p)\rightarrow \Phi (\frac{t_{0}\cdot t}{2\pi
},p).\right.
\end{eqnarray*}

Now, we pull-back the metric $g$ using $\Psi \;$or$\;\Phi $ and
obtain a metric on $\mathbf{R}\times L_{p}$ or
$\mathbf{S}^{1}\times L_{p}.$ Using \cite {Ponge}, it is easy to
see that this metric is $\varepsilon dt^{2}+g_{t},$ where
$g_{0}=g\mid _{L_{p}}.$
\end{proof}

\begin{remark}
\label{RemarkPropFacil}In the conditions of lemma
\ref{ProposFacil}, since the flow of $E$ takes leaves into leaves,
in order to ensure that all the integral curves with initial
condition on $L$ do not return to $L,$ is
sufficient to check this for only one of them, and therefore we obtain a $%
\mathbf{R}\times L$ type decomposition$.$

If we wish to obtain a $\mathbf{S}^{1}\times L$ type
decomposition$,$ we have to verify that all the integral curves
with initial condition on $L$ return to $L$ but intersect it at
only one point. But in this case, the existence of an integral
curve verifying the above property, does not guarantee it for the
others integral curves with initial condition on $L$.
\end{remark}

\begin{corollary}
\label{ColorPropFacil}Let $M$ be a semi-Riemannian manifold and
$E$ an unitary, irrotational and complete vector field. Take $p\in
M$ such that the integral curves of $E$ with initial value on
$L_{p}$ intersect $L_{p}$ at only one point.

\begin{itemize}
\item  If $E$ is orthogonally conformal then $M$ is isometric to
one of the
twisted product $\mathbf{R}\times _{f}L_{p}$ or $\mathbf{S}^{1}\times _{f}L_{p}$, $%
g=\varepsilon dt^{2}+f^{2}g_{0},$ where $g_{0}=g\mid _{L_{p}}$ and $%
f(t,x)=\exp (\int_{0}^{t}\frac{divE(\Phi _{x}(s))}{n-1}ds).$

\item  If $E$ is orthogonally conformal and $grad\,divE$ is proportional to $%
E,$ then $M\,$is isometric to one of the warped product
$\mathbf{R}\times _{f}L_{p}$ or $\mathbf{S}^{1}\times _{f}L_{p}$,
$g=\varepsilon
dt^{2}+f^{2}g_{0},$ where $g_{0}=g\mid _{L_{p}}$ and $f(t)=\exp (\int_{0}^{t}%
\frac{divE(\Phi _{p}(s))}{n-1}ds).$
\end{itemize}
\end{corollary}

\begin{proof}
Using the lemma \ref{ProposFacil}, $M$ is diffeomorphic to
$\mathbf{R}\times L_{p}$ or $\mathbf{S}^{1}\times L_{p}.$ If $E$
is orthogonally conformal, the orthogonal leaves are umbilic (see
lemma \ref{Propiedades de las dos Foliaciones}), and therefore we
obtain a twisted product $\mathbf{R}\times _{f}L_{p}$ or
$\mathbf{S}^{1}\times _{f}L_{p}$ with metric $\varepsilon
dt^{2}+f^{2}g_{0}$ where $g_{0}=g_{\left| L_{p}\right. }$ and
$f(0,q)\equiv 1.$ If we take $v\in U_{p}^{\perp },$ then
$\frac{divE}{n-1}\cdot v=\nabla _{v}E.$ But using the conexion
formulae of a twisted product \cite{Ponge} we get$\;\nabla
_{v}E=g(E,grad\log f)v,$ and so $\frac{divE}{n-1}=E(\log f).$ Thus
$f(t,x)=\exp (\int_{0}^{t}\frac{divE(\Phi _{x}(s))}{n-1}ds).$

If moreover $grad\,divE$ is proportional to $E,$ then $divE(\Phi
_{x}(s))=divE(\Phi _{p}(s))$ for all $x\in L_{p}$ and all $s\in
\mathbf{R}$.
Therefore $f(t)=\exp (\int_{0}^{t}\frac{divE(\Phi _{p}(s))}{n-1}ds).$%
\end{proof}

\begin{remark}
\label{RemarkLocalmente}Observe that the conclusion of lemma \ref
{ProposFacil} and corollary \ref{ColorPropFacil} are true locally
\cite{Montiel}. If $(a,b)\subset \mathbf{R}$, a warped product
$((a,b)\times _{f}L,-dt^{2}+f^{2}g_{L})$ is called a Generalized
Robertson-Walker spacetime \cite{Miguel}.
\end{remark}

\begin{example}
\label{EjemploHarris}If we suppose $M$ causal instead of the
condition about the integral curves with initial value on $L_{p}$
the conclusion of corollary \ref{ColorPropFacil} is not true,
compare with \cite{Harris2}.

We take $\tilde{M}=\mathbf{R}^{2}$ with the Minkowski metric and
the isometry $ \Phi (t,x)=(t,x+1)$. Let $\Gamma $ be the subgroup
of isometries generated by $\Phi $ and $M=\tilde{M}/ \Gamma $. We
consider $X=\sqrt{\frac{3}{2}}\frac{
\partial }{\partial t}+\sqrt{\frac{1}{2}}\frac{\partial }{\partial x}.$
Since $\Phi $ preserves the vector field $X,$ we can define the
vector field $U_{\Pi (p)}=\Pi _{*_{p}}(X_{p})$. Both $U$ and $X$
are parallel and complete. The manifold $M$ is causal, but it does
not split. This example can be trivially extended to any
dimension.
\end{example}

Simply connectedness implies that an irrotational vector field is
a gradient. If it has never null norm, it is immediate that the
integral curves meet the orthogonal leaves for only one value of
its parameter. So, with some additional hypothesis, we can use
lemma \ref{ProposFacil} to get a $(\mathbf{R}\times L,\varepsilon
dt^{2}+g_{t})$ type metric decomposition$.$ We can assume directly
that the vector field is a gradient and state the following: let
$M$ be a semi-Riemannian manifold and $f:M\rightarrow \mathbf{R}$
a function which gradient has never null norm and
$\frac{grad\,f}{\mid grad\,f\mid }$ is complete. If

\begin{itemize}
\item  $H^{_{f}}=0,$ then $M$ is isometric to a direct product
$\mathbf{R\times }L$.

\item  $H^{_{^{f}}}=a\cdot g,$ then $M$ is isometric to a warped product $%
\mathbf{R\times }L.$

\item  $H^{_{f}}=a\cdot g+bE^{*}\otimes E^{*},$ where $a,b\in
C^{\infty }(M), $ then $M$ is isometric to a twisted product
$\mathbf{R\times }L$.
\end{itemize}

But there are other ways to ensure that the integral curves with
initual values on a leaf does not return to the same leaf, as it
is shown in the following results.

\begin{corollary}
Let $M$ be a semi-Riemannian manifold and $U$ an irrotational and
conformal
vector field, with never null norm and complete unitary. Suppose that $%
\lambda $ is not constant. Then

\begin{enumerate}
\item  If $divU\geq 0$ (or $divU\leq 0)$ then $M$ is isometric to
a warped product$\;\mathbf{R\times }L.$

\item  If $Ric(U)\leq 0$ then $M$ is isometric to a warped product$\;\mathbf{%
R\times }L.$
\end{enumerate}
\end{corollary}

\begin{proof}
Since $\nabla U=a\cdot id,$ and $a=E(\lambda )$ it follows that
$divU=n\cdot E(\lambda )$ and $Ric(U)=-(n-1)U(E(\lambda )).$ Since
$\lambda $ is not constant, there is a point $p\in M$ such that
$divU_{p}\neq 0.$ Take $L$ the leaf through $p$ and $\gamma (t)$
an integral curve of $E$ with $\gamma (0)\in L.$

If $divU\geq 0$ (or $divU\leq 0),$ then $\lambda (\gamma (t))$ is
increasing (or decreasing), and since $\lambda $ is constant
through the leaves,$\gamma $ can not return to $L.$

If $Ric(U)\leq 0$ then $\frac{d}{dt^{2}}\lambda (\gamma (t))\geq 0,$ and so $%
\frac{divU_{p}}{n}t+\lambda (p)\leq \lambda (\gamma (t))$ for all $t\in \mathbf{%
R}$, and therefore $\gamma $ can not return to $L,$ since if this
happened then $\lambda (\gamma (t))$ would be periodic.

Therefore, in both cases the integral curves of $E$ with initial
condition
on $L$ intersect $L$ at only one value of its parameter. Since $divE=(n-1)%
\frac{E(\lambda )}{\lambda },$ it follows from corollary
\ref{ColorPropFacil} that, if we fix a point $p\in $ $M$, then $M$
is isometric to the warped
product $\mathbf{R\times }_{f}L$ where $f(t)=\frac{\lambda (\Phi _{p}(t))}{%
\lambda (p)}.$
\end{proof}
\begin{example}
Take $(\mathbf{R}\times _{e^{t}}N,-dt^{2}+e^{2t}g)$ where $(N,g)$
is a Riemannian manifold. We know that this warped product is not
timelike complete \cite{Miguel}. Then $U=e^{t}\frac{\partial
}{\partial t}$ is an irrotational and conformal vector field with
complete unitary and $divU>0.$ If $\Gamma $ is an isometry group
which preserves $U$ and the canonical foliations and the action is
properly discontinuous then $\left( \mathbf{R} \times
_{e^{t}}N\right) / \Gamma $ is a warped product manifold of $
\mathbf{R}\times _{e^{t}}L$ type.
\end{example}

\section{Irrotational vector fields with compact leaves}
\label{sec:3}

Completeness is a mild hypothesis but essential in most of
splitting theorems. We can give trivial counterexamples to these
theorems if we do not assume completeness. Some results change it
for the global hyperbolicity hypothesis \cite{Galloway}. We can
give the following theorems for irrotational vector fields with
compact leaves without assuming completeness.

\begin{theorem}
\label{TeoremaPrincipal}Let $M$ be a non compact semi-Riemannian
manifold and $E$ an unitary and irrotational vector field. Assume
$L_{p}$ is compact for all $p\in M.$ Then $M$ is isometric to
$(a,b)\times L,$ $-\infty \leq a<b\leq \infty ,$ where $E$ is
identified with $\frac{\partial }{\partial t}, $ $L$ is a compact
semi-Riemannian manifold and $g\;$is a parametrized
semi-Riemannian metric$\;\varepsilon dt^{2}+g_{t}.$
\end{theorem}
\begin{proof}
 Let $\Phi :A\subset \mathbf{R\times }M\rightarrow M$ be the
flow of $E.$ We know that $\Phi $ take any connected and open set
of a leaf into a leaf, see remark \ref{Piezas de Hojas a Hojas}.
Given $L$ a leaf we will show that the maximal definition interval
of $\Phi _{p}(t)$ is the same for all $p\in L.$ We know that for
each $p\in L$ there exists an open set $W_{p}\subset L$ and $\eta
_{p}$ with $(-\eta _{p},\eta _{p})\times W_{p}\subset A.$ Since
$L$ is compact, there is $\eta $ with $(-\eta ,\eta )\times
L\subset A.$

Let $(a,b)$ be the maximal interval such that $(a,b)\times
L\subset A.$ We claim that it is the maximal definition interval
of each integral curve with initial value on $L.$ In fact, suppose
that $\Phi _{t}(p_{0})$ is defined in
$(a,b+\delta )$ for some $p_{0}\in L.$ Then, there is a $\eta \ $such that $%
(-\eta ,\eta )\times L_{\Phi _{b}(p_{0})}\subset A.$ Since $\Phi _{-\frac{%
\eta }{2}}:L_{\Phi _{b}(p_{0})}\rightarrow L_{\Phi _{b-\frac{\eta
}{2}}(p)}$ is injective and a local diffeomorphism and $L_{\Phi
_{b}(p_{0})}$ is compact, it is a diffeomorphism. Therefore, given
$p\in L$ there is $q\in L_{\Phi _{b}(p_{0})}$ with $\Phi
_{-\frac{\eta }{2}}(q)=\Phi _{b-\frac{\eta }{2}}(p),$ and so $\Phi
_{t}(p)$ can be defined in $(a,b+\eta ).$ Since $ p\in L$ is an
arbitrary point, we obtain a contradiction$.$ Then the maximal
definition interval of the integral curves with initial values on
$L$ is $ (a,b).$ Now, we can define the local diffeomorphism

\begin{eqnarray*}
\Phi &:&(a,b)\times L\rightarrow M \\
&&\left. (t,p)\rightarrow \Phi _{t}(p).\right.
\end{eqnarray*}

We can show that $\Phi $ is onto as in lemma \ref{LemaFlujoSobre}.
Now we see that it is injective. It is sufficient to verify that
the integral
curves with intial values on $L$ does not meet $L$ again. If there is $%
(s,p)\in (a,b)\times L$ such that $\Phi _{s}(p)\in L,$ then
$M=\cup _{t\in [0,s]}\Phi _{t}(L),$ and $M$ would be compact.
Therefore, $\Phi $ is a diffeomorphism, and using \cite{Ponge}, we
can show that $\Phi ^{*}(g)=\varepsilon dt^{2}+g_{t},$where
$g_{0}=g\mid _{L}.$
\end{proof}

\begin{remark}
If $E$ is a complete vector field and a leaf is compact, then all
leaves are compact, since given two leaves $L_{p},L_{q}$ there is
a parameter $t\in \mathbf{R}$ such that $\Phi
_{t}:L_{p}\rightarrow L_{q}$ is a diffeomorphism. In this
situation $(a,b)=\mathbf{R}$.
\end{remark}

\begin{remark}
\label{remarkTeoremaHojasCompactas}In the same way as in collorary
\ref {ColorPropFacil}, if we assume $E$ unitary, irrotational and
orthogonally conformal in the above theorem, then $M$ is isometric
to a twisted product $ (a,b)\times L,$ and if moreover
$grad\,divE$ is proportional to $E$ we obtain a warped product
$(a,b)\times L.$

Observe that a chronological Lorentzian manifold is also non
compact \cite {Beem}.
\end{remark}

\begin{theorem}
Let $M$ be a semi-Riemannian manifold and $U$ an irrotational and
conformal vector field with never null norm such that $Ric(U)\geq
0.$ Assume that $ L_{p}$ is compact for all $p\in M$. Then $M$ is
isometric to a warped product $(a,b)\times L$ where $(a,b)\neq
\mathbf{R}$ and $E$ is identified with $\frac{\partial }{\partial
t}.$
\end{theorem}
\begin{proof}
We can prove in the same way that in the theorem
\ref{TeoremaPrincipal} that given an orthogonal leaf $L$, all the
integral curves with initial value on $ L$ have the same maximal
definition interval, say $(a,b).$ We also can show that $\Phi
:(a,b)\times L\rightarrow M$ is a local diffeomorphism which it is
into. Now, since $Ric(U)=-(n-1)U(E(\lambda ))$ it follows that
$(a,b)\neq \mathbf{R}.$ If $\Phi $ were not injective, since $\Phi
$ takes leaves into leaves, we would obtain that
$(a,b)=\mathbf{R}.$ So, $\Phi $ is a diffeomorphism and we can
show in the same way as in corollary \ref {ColorPropFacil} that
$M$ is isometric to the warped product $\mathbf{R}\times _{f}L$,
$g=\varepsilon dt^{2}+f^{2}g_{0},$ where $g_{0}=g\mid _{L}$ and $
f(t)=\frac{\lambda (\Phi _{p}(t))}{\lambda (p)}$ where $p$ is a
fixed point in $L.$
\end{proof}
Note that the Closed Friedmann Cosmological Model $(0,\pi )\times
_{f}\mathbf{S} ^{3}$ verifies the hypotesis of the above theorem
with the irrotational and conformal vector field
$U=f\frac{\partial }{\partial t}$.

\section{Application to Lorentzian Manifolds}
\label{sec:4}

We can use the above results to get decomposition theorems on
Lorentzian manifolds.

\begin{theorem}
Let $M$ be a Lorentzian manifold with positive sectional curvature
on timelike planes and $U$ a timelike, irrotational and conformal
vector field with complete unitary. Then $M$ is isometric to a
warped product $\mathbf{ R\times }L$ where $L$ is a Riemann
manifold and $E$ is identified with $ \frac{\partial }{\partial
t}.$
\end{theorem}
\begin{proof}
Take $L$ a leaf through and $p\in L.$ Given $v\in T_{p}L,$ a
direct computation gives us $K_{\{v,U_{p}\}}=\frac{E(E(\lambda
))}{\lambda }.$ Therefore, the sectional curvature of a plane
which contains $U_{p}$ only depends on $p.$ Let $\gamma $ be the
integral curve of $E$ with $\gamma (0)=p $ and take
$y:\mathbf{R\rightarrow }\mathbf{R}$ given by $y(t)=\lambda
(\gamma (t)).$ Then $K_{\{U_{\gamma (t)}\}}=\frac{y^{\prime \prime
}(t)}{y(t) },$ and if we define $f(t)=\log \frac{y(t)}{y(0)}$ we
obtain that $ 0<K_{\{U_{\gamma (t)}\}}=f^{\prime \prime
}+f^{\prime 2}.$ Now, it is easy to show that $f:[0,\infty
)\rightarrow \mathbf{R}$ has a finite number of zeros. If there
exists $t_{0}>0$ such that $\gamma (t_{0})\in L,$ then, since the
flow of $E$ takes leaves intos leaves, $\gamma (nt_{0})\in L$ for
all $n\in \mathbf{N}.$ But $\lambda $ is constant through the
leaves, thus $ \lambda (\gamma (nt_{0}))=\lambda (p)$ and
therefore $f(nt_{0})=0$ for all $ n\in \mathbf{N},$ which it is a
contradiction. Then, $\gamma $ does not return to $L,$ and using
corollary \ref{ColorPropFacil} we can conclude that $M$ is
isometric to the warped product $\mathbf{R}\times _{\frac{\lambda
(\gamma (t))}{ \lambda (p)}}L.$
\end{proof}

\begin{theorem}
Let $M$ be a Lorentzian manifold with positive sectional curvature
on all timelike planes and $U$ a timelike, irrotational and
conformal vector field with complete unitary. Then $M$ is
isometric to a warped product $\mathbf{R\times }L$ where $L $ is a
compact Riemann manifold and $E$ is identified with
$\frac{\partial }{
\partial t}.$
\end{theorem}

\begin{theorem}
Let $M$ be a complete and non compact Lorentzian manifold and $U$
a timelike, irrotational and conformal vector field. Suppose that
$Ric(v)\geq 0 $ for all $v\perp U$ and $\mid U\mid $ is not
constant and bounded from above. Then$\;M\;$is isometric to a
warped product $\mathbf{R\times }L$ where $L $ is a compact
Riemann manifold and $E$ is identified with $\frac{\partial }{
\partial t}.$
\end{theorem}
\begin{proof}
Observe that $E$ is complete, since it is geodesic. We show that
$X=\lambda ^{n}E$ is a complete vector field. Let us suppose it is
not true. Take $ \gamma :\mathbf{R\rightarrow }M$ an integral
curve of $E$ and $\beta :(c,d)\rightarrow M$ an integral curve of
$X$ with the same initial condition. Then $\beta (t)=\gamma
(h(t))$ where $h:(c,d)\rightarrow \mathbf{R}$ is a diffeomorphism.
But $d<\infty $ or $-\infty <c$ and $h^{\prime }(t)=\lambda
^{n}(\gamma (h(t)))$ is bounded$,$ which is a contradiction. Now,
we show that there is $q\in M$ with $\triangle \lambda (q)>0.$
Suppose $\triangle \lambda (q)\leq 0$ for all $q\in M.$ Since
$\nabla U=a\cdot id,$ where $a=E(\lambda ),\;\lambda $ is constant
through the leaves and $-\lambda ^{2}=g(U,U),$ we have
$grad\,\lambda =\frac{-a}{\lambda }U$ and
\begin{equation*}
\triangle \lambda =-E(a)-(n-1)\frac{a^{2}}{\lambda }\geq -E(a)-n
\frac{a^{2}}{\lambda }=-\frac{1}{\lambda ^{2n}}X(X(\lambda )).
\end{equation*}
Then $ X(X(\lambda ))\geq 0.$ Take $\beta :\mathbf{R\rightarrow
}M$ an integral curve of $X,$ and $y(t)=\lambda (\beta
(t)).\;$Then $0\leq y^{\prime \prime }$ and it is bounded from
above, but this is a contradiction. So, there is $q\in M$ with
$\triangle \lambda (q)>0.$ Consider $L_{q}$ the leaf through $q.$
Given a vector $v\in T_{z}L_{q}$ unitary for the metric $g$, a
direct computation gives us
\begin{equation*}
Ric_{L_{q}}(v)=Ric_{M}(v)-\frac{1}{\lambda }\left(
E(a)+(n-2)\frac{a^{2}}{\lambda }\right) \geq \frac{-1}{\lambda
}\left( E(a)+(n-1)\frac{a^{2}}{\lambda }\right) =\frac{\triangle
\lambda }{ \lambda }(z).
\end{equation*}
But since $\lambda $ is constant through
the leaf $L_{q},$ we deduce $Ric_{L_{q}}(v)\geq \frac{\triangle
\lambda }{\lambda }(q)>0.$ If we take the universal covering
$P:\tilde{M}\rightarrow M$ and $\tilde{U}$ the vector field such
that $P_{*_{e}}(\tilde{U}_{e})=U_{P(e)}$, then $\tilde{ U}$ is
irrotational and conformal too. Since $\tilde{M}$ is simply
connected, we know that it is isometric to a warped
product$\;\mathbf{R}\times \tilde{L}_{e},$ where $\tilde{L}_{e}$
is an orthogonal leaf of $\tilde{U} ^{\perp }.$ Since $M$ is
complete, $\tilde{M}$ is complete too, and so $ \tilde{L}_{e}$ is
complete \cite{Miguel}. If we take $e\in \tilde{M}$ such that
$P(e)=q$ then $P(\tilde{L}_{e})=L_{q}$ and since $P$ is a local
isometry $L_{q}$ is complete. Then, $L_{q}$ is a complete Riemann
manifold which satisfies $Ric_{L_{q}}(v)>c>0$ for all unitary
vector $v\in TL_{q}$ (for the induced metric on $L_{q}$) and so,
using Myers theorem \cite{Oneill} , it is compact. Since $\Phi
:\mathbf{R\times }L_{q}\rightarrow M$ is a local diffeomorphism
and it is onto, all the leaves are compact. Then, using theorem
\ref{TeoremaPrincipal}, $M$ is isometric to a warped product
$\mathbf{ R\times }L.$
\end{proof}
\begin{theorem}
\label{TeorGenerGR}Let $M$ be a non compact Lorentzian manifold
and $E$ a timelike, unitary and orthogonally irrotational vector
field such that the leaves of $E^{\perp }$ are compact and simply
connected. If $E(divE)\geq - \frac{(divE)^{2}}{n-1}$ and
$Ric(E)\geq 0$ then $M$ splits isometrically as a twisted product
$(a,b)\times _{f}L,$ where $L$ is an orthogonal leaf and $
f(t,p)=\frac{divE(p)}{n-1}t+1$, and the above inequalities are
equalities.
\end{theorem}
\begin{proof}
We take $p\in M$ and $\{e_{2},...,e_{n}\}$ an orthonormal basis of
$ E_{p}^{\perp }$ and we consider $A_{p}:E_{p}^{\perp }\rightarrow
E_{p}^{\perp }$ the endomorphism given by $A_{p}(X)=\nabla _{X}E.$
Since $E$ is orthogonally irrotational and timelike, $A_{p}$ is
diagonalizable. Thus $\frac{1}{n-1} tr(A_{p})^{2}\leq \mid \mid
A_{p}\mid \mid ^{2},$ where $tr(A_{p})$ denote the trace of
$A_{p}$ and $\mid \mid A_{p}\mid \mid
^{2}=\sum_{i=2}^{n}g(A_{p}(e_{i}),A_{p}(e_{i})),$ and the equality
holds if and only if $A_{p}(X)=\frac{tr(A_{p})}{n-1}X.$ Let
$\{E_{1},\ldots ,E_{n}\}$ be a frame near $p,$ with
$E_{1}(p)=E_{p}$ and $E_{i}(p)=e_{i}$. A straightforward
computation shows that
\begin{equation*}
Ric(E)_{p}=div\nabla _{E}E_{p}-E(divE)_{p}-\mid \mid A_{p}\mid
\mid ^{2}.
\end{equation*}

Using that $Ric(E)\geq 0$ and $E(divE)\geq
-\frac{(divE)^{2}}{n-1}$ we obtain that
\begin{equation*}
\mid \mid A_{p}\mid \mid ^{2}-\frac{(divE)^{2}}{n-1}\leq div\nabla
_{E}E.
\end{equation*}

But $divE_{p}=\sum_{i=2}^{n}g(\nabla _{e_{i}}E,e_{i})=tr(A_{p}).$
So,
\begin{equation*}
0\leq \mid \mid A_{p}\mid \mid ^{2}-\frac{1}{n-1}tr(A_{p})^{2}\leq
div\nabla _{E}E_{p}.
\end{equation*}

Since $p$ is arbitrary, $div\nabla _{E}E\geq 0$ on $M.$ Now, it is
known that through the leaves, the one form $g(\nabla _{E}E,\cdot
)$ is closed \cite{Walschap}. Let $L$ be the orthogonal leaf
trhough $p$. Since it is simply connected, there is a function
$f:L\rightarrow \mathbf{R}$ such that $ grad\,f=\nabla _{E}E.$
Now, a direct computation shows that
\begin{equation*}
\Delta _{L}e^{f}=e^{f}\cdot div\nabla _{E}E.
\end{equation*}

Since $L$ is compact and $\Delta _{L}e^{f}\geq 0$ on $L$, $f$ is
constant. Therefore $\nabla _{E_{p}}E=0$ and $\mid \mid A_{p}\mid
\mid ^{2}=\frac{1}{ n-1}tr(A_{p})^{2}$ but $p$ is arbitrary, thus
the above equalities remain valid on $M.$ So, $E$ is geodesic and
$\nabla _{X}E=\frac{tr(A)}{n-1}X$ for all $X\in E^{\perp }.$
Therefore, $E$ is unitary, irrotational and orthogonally
conformal. Now, using remark \ref{remarkTeoremaHojasCompactas},
$M$ is isometric to a twisted product $(a,b)\times _{f}L$ where
$f(t,p)=\exp (\int_{0}^{t}\frac{divE(\Phi p(s))}{n-1}ds).$ But,
since $\nabla _{E}E=0,$ the inequalities are equalities, so
$E(divE)=-\frac{(divE)^{2}}{n-1}$ and therefore
$f(t,p)=\frac{divE(p)}{n-1}t+1.$
\end{proof}
\begin{corollary}
Let $M$ be a non compact Lorentzian manifold and $E$ an unitary
and orthogonally irrotational vector field such that $Ric(E)\geq
0$ and $ E(divE)\geq 0.$ Assume that the orthogonal leaves are
compact and simply connected. Then $M$ is isometric to a direct
product $(a,b)\times L$, where $ L$ is a compact and simply
connected Riemann manifold.
\end{corollary}
\begin{proof}
Since $E(divE)\geq 0\geq -\frac{(divE)^{2}}{n-1}$ it follows from
the above theorem that $M$ is isometric to $(a,b)\times _{f}L,$
where $f(t,p)=\frac{ divE(p)}{n-1}t+1$, and the equality holds. It
is $0\leq E(divE)=-\frac{(divE)^{2}}{n-1}\leq 0.$ So $ divE=0$ and
$f(t,p)=1.$
\end{proof}
In a warped product $\mathbf{R}\times _{f}L$, if $\frac{\partial
}{\partial t}$ is complete and $Ric(\frac{\partial }{\partial
t})\geq 0$ then $f\equiv cte$ . On the other hand, in a twisted
product $\mathbf{R}\times _{f}L$, the same conditions on
$\frac{\partial }{\partial t}$ implies that $f$ is independent of
the variable $t$. Therefore, in the above theorem or corollary if
we assume that $E$ is complete, or $M$ timelike complete, then
$(a,b)=\mathbf{R}$ and we get a direct product. This shows that if
we want to get more general decomposition theorems with
$Ric(E)\geq 0$, then we must drop the completeness hypothesis.

A leaf is achronal if a timelike and future directed curve meets
the leaf at most once. Particularly, the integral curves of $E$
only meet the leaves one time. The achronality of the leaves is
more restrictive than the chronologicity$,$ and it is well known
that a chronological manifold is non compact. So, the achronality
of the leaves implies $M$ is non compact. Then, the above theorem
and corollary are generalizations of theorem 1 in \cite {Garcia3}.

Observe that the following twisted product verifies the hypothesis
of the theorem \ref{TeorGenerGR}. Take $(-1,\infty )\times
\mathbf{S}^{1}$ with the metric $g=-dt^{2}+f^{2}g_{can},$ where
the function is $f(t,e^{is})=t+2+\cos (s).$

\section{Irrotational vector fields with periodicity}
\label{sec:5}

Let $M$ be a semi-Riemannian manifold and $U$ an irrotational and
pregeodesic vector field with never null norm and complete
unitary. Then, the one form $w=g(\cdot ,U)$ is closed, because $U$
is irrotational. Following \cite{HectorHirsch}, we take the
homomorphism
\begin{eqnarray*}
\Psi &:&H_{1}(M,\mathbf{R})\rightarrow \mathbf{R} \\
&&\left. [\sigma ]\rightarrow \int_{\sigma }w.\right.
\end{eqnarray*}

Then, $G=\Psi (H_{1}(M,\mathbf{R}))$ is a subgroup of
$\mathbf{R},$ so there are three possibilities

\begin{enumerate}
\item  $G=0,$ and therefore $U$ is a gradient. Then $M$ is isometric to $(
\mathbf{R}\times L,g),$ where $g=\varepsilon dt^{2}+g_{t}.$

\item  $G\approx \mathbf{Z},$ and $M$ is a fibre bundle over
$\mathbf{S}^{1},$ with fibres the leaves of $U^{\perp }.$

\item  $G$ is dense in $\mathbf{R}.$
\end{enumerate}

Then, if for example, $\Pi _{1}(M)$ is finite, $U$ is a gradient
and $M=\mathbf{ R}\times L.$

The $\mathbf{S}^{1}\times L$ type decomposition is not frequent,
and it is more complicated than the $\mathbf{R}\times L$ type$,$
as it were commented in remark \ref{RemarkPropFacil}. The
following example shows the typical difficulty that presents this
type of decompositions.

\begin{example}
Take $\mathbf{R}\times _{f}\mathbf{S}^{3}(\frac{1}{2}),$
$g=dt^{2}+f^{2}g_{can},$ where $f(t)=\sqrt{3+\sin (2t)}.$ Since
each factor is Riemannian and complete, this warped product is
complete. Take the isometry $\eta :\mathbf{R} \times
_{f}\mathbf{S}^{3}(\frac{1}{2})\rightarrow \mathbf{R}\times
_{f}\mathbf{S}^{3}( \frac{1}{2})$ given by $\eta (t,p)=(t+\pi
,-p).$ We call $\Gamma $ the isometry group generated by $\eta .$
Then it is easy to check that $\Gamma $ acts in a properly
discontinuous manner. We consider the quotient $\Pi :\mathbf{
R}\times _{f}\mathbf{S}^{3}(\frac{1}{2})\rightarrow M=\left(
\mathbf{R}\times _{f} \mathbf{S}^{3}(\frac{1}{2})\right) / \Gamma
$ and take the irrotational and conformal vector field
$V=f\frac{\partial }{\partial t}.$ Since $\eta $ preserves the
vector field $V,$ there is a vector field $U$ on $M$ such that
$\Pi _{*}(V)=U$ and it is irrotational and conformal too. The
integral curves of $U$ are periodic, but $M$ is not isometric to a
product $\mathbf{S} ^{1}\times L$ since the integral curves of $U$
intersect each leaf in two different points. In fact, take $\Pi
(0,p)=\Pi (\pi ,-p)\in M$ and call $ L=\Pi \left( \{0\}\times
\mathbf{S}^{3}(\frac{1}{2})\right) =\Pi \left( \{\pi \}\times
\mathbf{S}^{3}(\frac{1}{2})\right) $ the leaf through it. Observe
that $\Pi (t,p)$ is the integral curve of $\frac{U}{\left|
U\right| }$ through $ \Pi (0,p)$ and it intersects the above leaf
at $\Pi (0,p)$ and $\Pi (\pi ,p)$ .
\end{example}

A foliation is regular if for each $p\in M$ there exists an
adapted coordinated system such that each slice belongs to a
unique leaf \cite {Palais}.

\begin{theorem}
Let $M$ be a chronological Lorentzian manifold and $U$ a timelike,
irrotational and conformal vector field with complete unitary.
Suppose that the foliation $U^{\perp }$ is regular and let $L$ be
a leaf of $U^{\perp }$. Then $M$ is isometric to a warped product
$\mathbf{R\times }L$ or there is a lorentzian covering map $\Psi
:M\rightarrow \mathbf{S}^{1}\times N,$ where $N$ is a quotient of
$L$ and $\mathbf{S}^{1}\times N$ is a warped product.
\end{theorem}
\begin{proof}
If the integral curves of $E$ with initial value on $L$ do not
meet $L$ again, we know that $M$ is isometric to $\mathbf{R\times
}_{\frac{\lambda (\Phi _{p}(t))}{\lambda (p)}}L$ with $p\in L$
(corollary \ref{ColorPropFacil}).

Suppose there is an integral curve $\gamma $ that meets $L$ again.
Now, we can define $ t_{0}=\inf \{t>0:\gamma (t)\in L\}.$ Since
$U^{\perp }$ is a regular foliation, $t_{0}>0$ and it is a
minimum. Since $\Phi $ takes leaves into leaves it is easy to
verify that $\Phi _{t_{0}}(q)\in L_{q}$ for all $q\in M$ . Now,
since $U$ is conformal, the diffeomorphisms $\Phi
_{t}:L_{q}\rightarrow L_{\Phi _{t}(q)}$ are conformal with
constant factor $ \left( \frac{\lambda (\Phi _{t}(q))}{\lambda
(q)}\right) ^{2}.$ But $\Phi _{t_{0}}(q)\in L_{q}$ and $\lambda $
is constant through the leaves, so $ \Phi
_{t_{0}}:L_{q}\rightarrow L_{q}$ is an isometry. Since $\Phi
_{t_{0}}$ preserves the vector field $E,$ we have that $\Phi
_{t_{0}}:M\rightarrow M$ is an isometry.

Let $\Gamma $ be the subgroup of isometries generates by $\Phi
_{t_{0}}$. We can suppose that $U$ is future pointing. Since $M$
is chronological, $\Gamma $ is isomorphic to $\mathbf{Z}.$

Now we show that given $q\in M$ there is an open set $B$ with
$q\in B,$ such that for all $z\in B$ the integral curve of $E$
with initial value $z$ leaves $B$ before $t_{0}$ and it does not
return to $B$.

We know that there is an open set $B$ such that $\Phi
:(-\varepsilon ,\varepsilon )\times _{\frac{\lambda (\Phi
_{q}(t))}{\lambda (q)} }W\rightarrow B$ is an isometry, where
$W\subset L_{q}$ (see remark \ref {RemarkLocalmente}). We can
assume that $W$ is the convex ball $B_{q}(\delta )$ in $L_{q}$
with $\delta <\frac{\varepsilon ^{2}}{4}.$ Suppose that there is
$z\in B$ such that the integral curve $\Phi _{t}(z)$ returns to
$B.$ Then, since $E$ is indentified with $\frac{\partial
}{\partial t},$ there are $a,b\in W$ such that $\Phi _{s}(a)=b$
for some $s\in \mathbf{R}.$ Take $ \alpha :[0,1]\rightarrow W$ a
geodesic in $W$ with $\alpha (0)=b$ and $ \alpha (1)=a.$ Now we
consider the curve $\beta (t)=\Phi (-\frac{\varepsilon
}{2}(1-t),\alpha (t)),$ $t\in [0,1].$ This curve joins $\Phi
(-\frac{ \varepsilon }{2},b)$ with $\Phi (0,a)=a,$ and
\begin{equation*}
g(\beta ^{\prime }(t),\beta ^{\prime }(t))=-\frac{\varepsilon
^{2}}{4}+g(\alpha ^{\prime }(t),\alpha ^{\prime
}(t))<-\frac{\varepsilon ^{2}}{4}+\delta <0.
\end{equation*}
So, using the curve $ \Phi _{t}(a)$ and $\beta (t)$ we can
construct a piecewise smooth closed timelike curve, which is a
contradiction with the chronological hipothesys.

Now, we take the action of $\Gamma .$ on $M$. Let us see that this
is a properly discontinuous action. We have to see:

\begin{enumerate}
\item  Given $p\in M$ there exists an open set $U$, $p\in U,$ such that $
U\cap \Phi _{nt_{0}}(U)=\emptyset $ for all $n\in \mathbf{Z}.$

It is sufficient to take an open set B with the above property.
\item  Given $p,q\in M$ with $p\neq \Phi _{nt_{0}}(q)$ for all
$n\in \mathbf{Z}, $ there is open sets $U,V$ such that $p\in U,$
$q\in V$ and $U\cap \Phi _{nt_{0}}(V)=\emptyset \;$for all$\;n\in
\mathbf{Z}.$
\end{enumerate}

We suppose that this is not true and show that there is $m_{n}\in
\mathbf{N}$ such that $\left. \lim_{n\rightarrow \infty }\Phi
_{m_{n}t_{0}}(q)=p.\right. $

We take $U_{n}=\Phi ((-\frac{1}{n},\frac{1}{n})\times
B_{p}(\frac{1}{2n}))$ and $V_{n}=\Phi
((-\frac{1}{n},\frac{1}{n})\times B_{q}(\frac{1}{2n})).$ Since
property $2$ is not true, there is $m_{n}\;$with $U_{n}\cap \Phi
_{m_{n}t_{0}}(V_{n})\neq \emptyset .$ Using the fact that $\Phi
_{m_{n}t_{0}}:L_{q}\rightarrow L_{q}$ are isometries, it is easy
to verify that
\begin{equation*}
 \Phi _{m_{n}t}(V_{n})=\Phi
((-\frac{1}{n},\frac{1}{n})\times B_{\Phi
_{m_{n}t}(q)}(\frac{1}{2n})),
\end{equation*}
then it follows that $\lim_{n\rightarrow \infty }\Phi
_{m_{n}t_{0}}(q)=p.$

We claim that $m_{n}$ is constant from a $n_{1}$ forward$.$ If
this were not true, we take the open set $B$ with $p\in B$, such
that the integral curves of $E$ with initial values in $B,$ leave
it before $t_{0}$ and it does not return to it. Since $p\in B$,
there is $n_{0}$ such that if $n\geq n_{0}$ we have $\Phi
_{m_{n}t_{0}}(q)\in B.$ But, there are $m_{r},m_{s}$ such that $
m_{r}-m_{s}=k\geq 1$ and $m_{r},m_{s}\geq n_{0}$, thus $\Phi
_{kt_{0}}(\Phi _{m_{s}t_{0}}(q))$ is outside $B$ and this is in
contradiction with $\Phi _{kt_{0}}(\Phi _{m_{s}t_{0}}(q))=\Phi
_{m_{r}t_{0}}(q)\in B.$

Therefore $\Phi _{kt_{0}}(q)=p$ for some $k,$ and this is a
contradiction.

Now we take the quotient $P:M\rightarrow M/ \Gamma .$ We can take
a metric on $M/ \Gamma $ such that $P$ is a local isometry. Since
$ \lambda \;$is constant through the orthogonal leaves, $\Phi
_{t_{0}}$ preserves the vector field $U.$ So there is a timelike,
irrotational and conformal vector field $Y$ on $M/ \Gamma $ such
that $ P_{*_{e}}(U_{e})=Y_{p(e)}.$ The integral curves of $Y$
intersect the leaf of $Y^{\perp }$ given by $p(L)=N$ at only one
point, and since the integral curves of $Y$ are diffeomorphic to
$\mathbf{S}^{1},$ $M/ \Gamma \;$is isometric
to$\;\mathbf{S}^{1}\times \frac{\lambda }{\lambda (p)}N.$ It is
easy to prove that $\Gamma $ acts on $L$ and $L / \Gamma =N.$
\end{proof}

\begin{remark}
In the above theorem, we can not expect that the covering map
would be a diffeomorphism because $\mathbf{S}^{1}\times N$ is not
choronological. On the other hand, the example \ref{EjemploHarris}
satisfies the conditions of the above theorem, and we obtain the
covering map $p:M\rightarrow \mathbf{T}^{2}.$
\end{remark}

Given a foliation of arbitrary dimension we can define the
holonomy of a leaf of the foliation \cite{Camacho}. In some sense
it measures how intertwine the leaves through a small transversal
manifold around a fixed point. If the foliation is defined by the
integral curves of a vector field, then the holonomy is given by
its flow. In this case, an integral curve without holonomy means
that it is diffeomorphic to $\mathbf{R}$ or if it is periodic,
then all other nearby integral curves are periodic. Using this
notion, we can prove the following result.

\begin{theorem}
Let $M$ be a compact and orientable Riemann manifold with odd
dimension and let $U$ be an irrotational and conformal vector
field. Suppose that $K_{\Pi }\geq 0$ for all planes $\Pi \perp U,$
the norm $\mid U\mid $ is not constant and the integral curves are
without holonomy. Then $M$ is isometric to a warped product
$\mathbf{S}^{1}\times _{f}L$ where $L$ is a compact and simply
connected Riemannian manifold$.$
\end{theorem}
\begin{proof}
We know that $E$ is geodesic and therefore complete. If we take
the universal covering $\tilde{M}$ we know that it is isometric to
$\mathbf{R} \times _{\frac{\lambda }{\lambda (e)}}\tilde{L}_{e}.$
Since $M$ is compact it follows that $ \mathbf{R}\times
_{\frac{\lambda }{\lambda (e)}}\tilde{L}_{e}$ is complete. Then
$\tilde{L}_{e}$ is complete \cite{Miguel}, and so is $L.$ Now we
show that there exists a compact leaf with positive sectional
curvature. Let $ \gamma $ be an integral curve of $E.$ If $\gamma
$ does not meet the leaf $ L_{\gamma (0)}$ again then the integral
curves with initial condition on $ L_{\gamma (0)}$ do not meet
again $L_{\gamma (0)}$ and $M$ woud be diffeomorphic to
$\mathbf{R}\times L_{\gamma (0)}$ which is a contradiction with
the compacity of $M$. Then $\gamma $ meets again $L_{\gamma (0)}$
and therefore $f(t)=\lambda (\gamma (t))$ is periodic and non
constant. Then there exists $s$ with $f^{\prime }(s)>0.$ Now, let
$L$ be the leaf trough $ \gamma (s).$ If $\Pi $ is a plane of $L$
then we obtain $K_{\Pi }^{L}=K_{\Pi }^{M}+(\frac{E(\lambda
)}{\lambda })^{2}\geq (\frac{f^{\prime }(s)}{f(s)} )^{2}>0.$ Since
$L$ is a complete and orientable Riemann manifold with even
dimension and $K_{\Pi }^{L}>c>0$ for all planes $\Pi $, it follows
that it is compact and simply connected \cite{DoCarmo}. Take $p\in
L$ and consider $t_{0}=\inf \{t>0:\Phi _{t}(p)\in L\}.$ Now we
show that $t_{0}>0.$ Suppose $t_{0}=0.$ Then there exists
$t_{n}\rightarrow 0,$ $t_{n}>0,$ such that $p_{n}=\Phi
(t_{n},p)\in L.$ Since $L$ is compact we can assume that $p_{n}$
converges, necessarily to $p.$ We know that $\Phi :(-\varepsilon
,\varepsilon )\times W\rightarrow \theta $ is a diffeomorphism,
where $W$ is an open set in $L.$ So, we can suppose that $
p_{n}\in W,$ but then $t_{n}=0$ for all $n,$ and this is a
contradiction. Now it is easy to verify that $t_{0}$ is a minimum
and it is the minimun value which $\Phi _{t}(q)\in L$ for all
$q\in L.$ Since $U$ is irrotational and conformal, $\Phi _{t}$ are
conformal diffeomorphism with constant factor $\left(
\frac{\lambda (\Phi _{t}(p))}{\lambda (p)}\right) ^{2}.$ Then,
$\Phi _{t_{0}}:L\rightarrow L$ is an isometry. If $\omega $ is the
volume form of $ M$ then $i_{U}\omega $ is a volume form of $L.$
It is easy to verify that $ \Phi _{t_{0}}$ preserve this
orientation. Now, using the theorem of Synge \cite{DoCarmo}, we
can ensure that there exists $q\in L$ such that $\Phi
_{t_{0}}(q)=q.$ Since the integral curves have not holonomy, $\Phi
_{t_{0}}$ is the indentity near $q,$ but since it is an isometry,
$\Phi _{t_{0}}=id.$ Since the integral curves with initial
condition on $L$ intersect it at only one point, it follows from
\ref{ColorPropFacil} that $M$ is isometric to the warped product
$\mathbf{S}^{1}\times _{\frac{\lambda (\Phi _{p}(t))}{\lambda (p)
}}L$ where $p\in L.$
\end{proof}

\textbf{Acknowledgements} The authors acknowledge to Prof. S. G.
Harris for his valuable comments and suggestion concerning example
\ref{EjemploHarris}. The first author has been partially supported
by DGICYT Grant BFM2001-1825.

\end{document}